\documentclass[12pt,a4paper,twoside]{article}
\usepackage{amssymb}
\usepackage{amsmath}
\usepackage[mathscr]{eucal}

\date{}
\newcommand{\A}[1]{\vspace{3mm}}
\newcommand{\D}{{\rm d}}
\newcommand{\E}{{\rm e}}
\newcommand{\I}{{\rm i}}
\newcommand{\R}{\mathbb{R}}
\begin{document}

\author{Rolf Schneider and Franz E. Schuster}
\title{Rotation Invariant Minkowski Classes\\of Convex Bodies}

\maketitle

\begin{abstract}
\noindent A Minkowski class is a closed subset of the space
of convex bodies in Euclidean space $\R^n$ which is closed
under Minkowski addition and non-negative dilatations. A convex body
in $\R^n$ is universal if the expansion of its support function
in spherical harmonics contains non-zero harmonics of all orders.
If $K$ is universal, then a dense class of convex bodies $M$ has the following property.
There exist convex bodies $T_1,T_2$ such that $M+T_1=T_2$, and $T_1,T_2$ belong to
the rotation invariant Minkowski class generated by $K$.
We show that every convex body $K$ which is not centrally
symmetric has a linear image, arbitrarily close to $K$, which is universal.
A modified version of the result holds for centrally symmetric convex bodies.
In this way, we strengthen a result of S. Alesker, and at the same time
give a more elementary proof.\\[2mm]
\end{abstract}

\renewcommand{\thefootnote}{{}}
\footnote{This work was supported by the European Network PHD,
FP6 Marie Curie Actions, RTN, \linebreak Contract
MCRN-2004-511953. The second author was also supported by the
Austrian Science Fund (FWF), within the project ``Affinely
associated bodies'', Project Number P16547-N12. }

\section{Introduction and Main Results}

\noindent Let ${\cal K}^n$ denote the space of convex bodies
(non-empty, compact, convex sets) in \linebreak[4] $n$-dimensional
Euclidean space ${\mathbb R}^n$ ($n\ge 2$). The basic
algebraic structures on $\mathcal{K}^n$ are Minkowski
addition and dilatation, defined by
\[K+L:=\{x+y:x \in K, \;y \in L\} \qquad \mbox{and} \qquad \lambda K:= \{\lambda x:x \in K\},\]
respectively, for $K,L\in{\cal K}^n$ and $\lambda\ge 0$.
By a {\em Minkowski class} in $\R^n$ we understand (slightly modifying
the definition of an $M$-class given in \cite[p. 164]{Sch93})
a subset of ${\cal K}^n$ which is closed in the Hausdorff metric and closed under
Minkowski addition and dilatation.

If $G$ is a group of transformations of ${\mathbb R}^n$, the
Minkowski class ${\cal M}$ is called $G$-{\em invariant} if
$K\in{\cal M}$ implies $gK\in{\cal M}$ for all $g\in G$. The
smallest $G$-invariant Minkowski class containing a given convex
body $K\in{\cal K}^n$ is said to be the $G$-invariant Minkowski
class {\em generated by $K$}. It consists of all convex bodies
which can be approximated
by bodies of the form $\lambda_1g_1K+\ldots+\lambda_kg_kK$
with $k\in {\mathbb N}$, $\lambda_1,\dots,\lambda_k\ge0$,
and $g_1,\dots, g_k\in G$.

The elements of a Minkowski class ${\cal M}$ will be called ${\cal M}$-{\em bodies}.
Further, the convex body $K$ is a {\em generalized ${\cal M}$-body}
if there exist ${\cal M}$-bodies $T_1,T_2$ such that $K+T_1=T_2$.

Let ${\cal M}$ be the $G$-invariant Minkowski class
generated by a convex body $K$. Then the convex body $M$ is a generalized ${\cal M}$-body
if and only if its support function can be approximated, uniformly
on the unit sphere, by functions from the vector space spanned by
the support functions of the $G$-images of $K$.

We recall a classical example. A convex body in ${\mathbb R}^n$
($n\ge 2$) is a {\em zonoid} if it can be approximated by
Minkowski sums of finitely many closed line segments. Thus, the set
${\cal Z}^n$ of zonoids is the rigid motion invariant Minkowski
class generated by a (non-degenerate) segment. Since the
affine image of a segment is a segment, ${\cal Z}^n$ is also the
affine invariant Minkowski class generated by a segment.

Every zonoid belongs to the subset ${\cal K}^n_c\subset{\cal K}^n$ of convex bodies
which have a centre of symmetry; such bodies are called {\em symmetric}
in the following, and {\em origin symmetric} if $0$ is the centre of symmetry. It is easy to see that ${\cal
Z}^2={\cal K}_c^2$, but for $n\ge 3$ the set ${\cal Z}^n$ is
nowhere dense in ${\cal K}_c^n$. A convex body $K\subset{\mathbb
R}^n$ is called a {\em generalized zonoid} if there exist zonoids
$Z_1,Z_2$ such that $K+Z_1=Z_2$. The set of generalized zonoids turns
out to be dense in ${\cal K}_c^n$, see \cite[Corollary 3.5.6]{Sch93}.
Generalized zonoids played a critical role in the first author's
solution \cite{Sch67} of the Shephard problem and in Klain's
classification of translation invariant, even and simple
valuations, see \cite{Kla95}, \cite{Kla99}. More information on
zonoids and generalized zonoids is found in the survey articles
\cite{SW83}, \cite{GW93} and in Section 3.5 of the book
\cite{Sch93}.

In the following, we want to replace the segment, which is used in
the definition of zonoids, by other convex bodies. The
non-denseness of zonoids mentioned above extends as follows. For
$n\ge 3$, the affine invariant Minkowski class generated by a
convex body (a symmetric convex body) is nowhere dense
in ${\cal K}^n$ (nowhere dense in ${\cal K}_c^n$). This follows
from \cite[Theorem 3.3.3]{Sch93}.

Our main issue here is the question analogous to the denseness of
generalized zonoids. Now we have to distinguish between convex
bodies with or without a centre of symmetry. Let $\mathcal{K}_o^n
\subset \mathcal{K}_c^n$ be the subset of origin symmetric convex
bodies. In the following, a convex body is called {\em non-symmetric}
if it does not have a centre of symmetry, and
{\em non-trivial} if it has more than one point.
Alesker \cite{Ale03} has proved the following theorem, in a
different but equivalent formulation.

\A

\noindent{\bf Theorem (Alesker).} (a) {\em If ${\cal M}$ is the
$SL(n)$-invariant Minkowski class generated by a non-symmetric
convex body, then the set of generalized ${\cal M}$-bodies is
dense in
${\cal K}^n$.}\\[2mm]
(b) {\em Let ${\cal M}$ be the $SL(n)$-invariant Minkowski class
generated by a non-trivial symmetric convex body $K$ with centre
different from $0$ $($with centre $0$$)$. Then the set of generalized
${\cal M}$-bodies is dense in ${\cal K}_c^n$ $($dense in ${\cal K}^n_o$$)$. }

\A

Part (b) of this theorem extends the statement about generalized
zonoids recalled above (and can be deduced from it, see
\cite[Remark 9]{Ale03}).

Note that in Alesker's result in effect the general linear group $GL(n)$ is
applied to the convex body $K$, since Minkowski classes are dilatation invariant. In contrast to the case of segments, for
a general convex body $K$ the affine invariant Minkowski class
does not coincide with the rigid motion invariant
Minkowski class generated by $K$. However, part (a) of Alesker's theorem
should be compared to an immediate consequence of a result proved and used in \cite{Sch96}:

\A

\noindent{\bf Theorem.} {\em Let $T\subset {\mathbb R}^n$ be a
triangle. Then there exists an affine map $A$ such that for the
rotation invariant Minkowski class ${\cal M}$ generated by $AT$, the
set of generalized ${\cal M}$-bodies is dense in ${\cal K}^n$.}

\A

Hence, for an arbitrary triangle $T$, the $SL(n)$-invariant Minkowski class
generated by $T$ in part (a) of Alesker's Theorem may be replaced by the rotation invariant Minkowski class
generated by $AT$, for one suitably chosen affine map $A$. Alesker \cite{Ale03}, p. 58,
remarks that it is not clear whether this result can be obtained
by his method.

Our aim in the following is to strengthen Alesker's theorem in a
way suggested by the latter theorem. Instead of applying
all linear transformations to a given convex body, it is sufficient to perturb it only a
little by an appropriate linear map and then to apply only rotations and dilatations.

\A

\noindent{\bf Theorem 1.} (a) {\em Let $K\in{\cal K}^n$ be a
non-symmetric convex body. Then there exists a linear map $A$,
arbitrarily close to the identity, such that for the rotation
invariant Minkowski class ${\cal M}$ generated by $AK$, the set of
generalized ${\cal M}$-bodies is dense in} ${\cal K}^n$.\\[2mm]
(b) {\em Let $K\in{\cal K}^n$ be a non-trivial symmetric convex body, with centre
different from $0$ $($with centre $0$$)$.
Then there exists a linear map $A$, arbitrarily close to the identity, such that for the rotation
invariant Minkowski class ${\cal M}$ generated by $AK$, the set of
generalized ${\cal M}$-bodies is dense in ${\cal K}_c^n$ $($dense in ${\cal K}_o^n$$)$}.

\A

Clearly, the perturbation by the linear transformation $A$ is necessary
in general, as shown by the case of a ball in the symmetric case, and by a
body of constant width in the non-symmetric case.

Besides strengthening Alesker's theorem, our second aim was to
give an easier proof for it. Whereas \cite{Ale03} employs deep
results on representations of the general linear group (proving
an irreducibility theorem which is analogous to Alesker's \cite{Ale01}
irreducibility theorem used in the theory of valuations), our proof
uses spherical harmonics and is comparatively elementary and self
contained.

The basic notion in our method of proof is that of
universal convex bodies. A convex body $K$ in ${\mathbb R}^n$ is
called {\em universal} ({\em centrally universal}\,) if the expansion of its support function
in spherical harmonics contains non-zero harmonics of all orders (of all even orders).
Universal convex bodies were introduced and applied by the first
author in \cite[p. 70]{Sch74}. The following theorem shows why
they are crucial for our result.

\A

\noindent {\bf Theorem 2.} {\em Let $K\in{\cal K}^n$ be a convex body, and let ${\cal M}$ be the rotation invariant Minkowski class
generated by $K$.}\\[2mm]
(a) {\em The set of generalized ${\cal M}$-bodies is dense in
${\cal K}^n$ if and only if $K$ is universal. }\\[2mm]
(b) {\em Let $K$ be symmetric. If $K$ has centre different from $0$
$($centre $0$$)$, then the set of generalized ${\cal M}$-bodies
is dense in ${\cal K}_c^n$ $($dense in ${\cal K}_o^n$$)$ if and only
if $K$ is centrally universal. }

\A

The only explicit convex bodies which are known to be universal
are the triangles with Steiner point not at the origin
and with the property that at least one of their angles is an
irrational multiple of $\pi$, see \cite{Sch74}, p. 71.

We consider the next theorem as our main result. In view of Theorem 2, it
implies Theorem 1.

\A

\noindent{\bf Theorem 3.} (a) {\em Let $K\in{\cal K}^n$ be a non-symmetric convex
body. Then there exists a linear
transformation $A$, arbitrarily close to the identity, such that
$A K$ is universal.}\\[2mm]
(b) {\em Let $K\in{\cal K}^n$ be a non-trivial convex body. Then there exists a linear
transformation $A$, arbitrarily close to the identity, such that
$A K$ is centrally universal.}

\A

\noindent{\bf Remark.} Our proof of Theorem 3, and thus of Theorem 1,
shows that, in fact, in every neighbourhood of the identity in $GL(n)$, almost all
(in the sense of measure) linear maps have the required property.

We will prove Theorem 2 in Section 2, part (b) of Theorem 3 in Section 3, and
finish the proof of Theorem 3 in Section 4.

\section{Universal Convex Bodies}

In this section, we collect a few facts about spherical harmonics
and convex bodies. An introduction to spherical harmonics and their use in
convexity is found in the book of Groemer \cite{Gro96}; see also
the short appendix of \cite{Sch93}.

By $S^{n-1}$ we denote the unit sphere of $\R^n$ and by $\sigma$
the spherical Lebesgue measure on $S^{n-1}$. A spherical harmonic
of dimension $n$ and order $m$ is the restriction to $S^{n-1}$ of
a harmonic polynomial of degree $m$ on $\R^n$.
For $m\in{\mathbb N}_0$, we denote by ${\cal H}^n_m$ the real
vector space of spherical harmonics of dimension
$n$ and order $m$. $\mathcal{H}^n$ will denote the space of all
finite sums of spherical harmonics of dimension $n$.

${\cal H}^n_m$ is a finite-dimensional subspace of $C(S^{n-1})$,
the vector space of real continuous functions on
$S^{n-1}$; let $N(n,m)$ denote its dimension. With respect to the
scalar product defined by
\[ (f,g) :=\int_{S^{n-1}} fg\,\D\sigma,\qquad f,g\in C(S^{n-1}),\]
spherical harmonics of different orders are
orthogonal. By
\[ \pi_m:C(S^{n-1})\to {\cal H}^n_m \]
we denote the orthogonal projection to ${\cal H}^n_m$. In each
space ${\cal H}^n_m$ we choose an orthonormal basis
$\{Y_{m1},\dots,Y_{mN(n,m)}\}$, which will be kept fixed in the
following; then
\begin{equation}\label{2.2}
\pi_m f= \sum_{j=1}^{N(n,m)} (f,Y_{mj})Y_{mj} \qquad\mbox{for }f\in C(S^{n-1}).
\end{equation}
One also writes
\[ f \sim \sum_{m=0}^{\infty} \pi_m f \]
and calls this the {\em condensed harmonic expansion} of $f$
(Groemer \cite{Gro96}, p. 72). The series converges to $f$ in the
$L_2$-norm.

The rotation group $SO(n)$ acts on $C(S^{n-1})$ by means of
$(\vartheta f)(u)=f(\vartheta^{-1}u)$, $u\in S^{n-1}$.
The space $\mathcal{H}^n_m$ is invariant under rotations. Thus,
for any rotation $\vartheta \in SO(n)$, we have
\begin{equation} \label{tij}
\vartheta Y_{mj}(u) = \sum
\limits_{i=1}^{N(n,m)}t_{ij}^m(\vartheta)Y_{mi}(u), \qquad u \in
S^{n-1},
\end{equation}
with real coefficients $t_{ij}^m(\vartheta)$. Let $\nu$ denote
the normalized Haar measure on the compact group $SO(n)$. The
following formula was proved in \cite[Lemma 3]{Sch96}. If $f \in
C(S^{n-1})$, then
\begin{equation} \label{form1}
\int_{SO(n)}\vartheta f(u)t_{ij}^m(\vartheta)\,\D\nu(\vartheta) =
N(n,m)^{-1} (f,Y_{mj}) Y_{mi}(u)
\end{equation}
for $u \in S^{n-1}, n \in \mathbb{N}_0$, and $i,j = 1, \ldots,
N(n,m)$.

A convex body $K \in \mathcal{K}^n$ is determined by its support
function $h(K,\cdot)$, defined on $\mathbb{R}^n$ by
$h(K,x)=\max\{ \langle x,y \rangle :y \in K\}$. Since $h(K,\cdot)$ is positively
homogeneous of degree one, it is determined by its restriction
to the sphere $S^{n-1}$, which we
denote by $h_K$. If $K, L \in \mathcal{K}^n$, then $h(K +
L,\cdot)=h(K,\cdot) + h(L,\cdot).$ Moreover, convergence in the
Hausdorff metric on $\mathcal{K}^n$ is equivalent to uniform
convergence of support functions on $S^{n-1}$.

The functions $\pi_m h_K$, $m\in{\mathbb N}$, determine the convex
body $K$ uniquely. In particular, a convex body $K$ is
symmetric if and only if $\pi_m h_K=0$ for all odd numbers
$m\not=1$, and $K$ is one-pointed if and only if $h_K \in
\mathcal{H}_1^n$. It follows that for $m\not=1$, the projection
$\pi_m h_K$ is invariant under translations of $K$. We also note
two special cases:
\begin{equation*}
 (\pi_0 h_K)(u) = \frac{1}{2}b(K) \qquad \mbox{and} \qquad
(\pi_1 h_K)(u)=\langle s(K),u\rangle = h_{\{s(K)\}}(u)
\end{equation*}
for $u \in S^{n-1}$ and $K \in \mathcal{K}^n$. Here, $b:
\mathcal{K}^n \rightarrow \mathbb{R}$ is the mean width of the
convex body $K$, defined by
\[b(K)=\frac{2}{\omega_n} \int_{S^{n-1}}h_K\,\D \sigma, \]
where $\omega_n=\sigma(S^{n-1})$. The rigid motion equivariant map $s: \mathcal{K}^n
\rightarrow \mathbb{R}^n$ is the Steiner point map, defined by
\[s(K)=\frac{n}{\omega_n}\int_{S^{n-1}} h_K(u)u\,\D\sigma(u). \]
It follows that $\pi_0 h_K=0$ if and only if $K$ contains only one
point, and $\pi_1 h_K=0$ if and only if $s(K)=0$.

We can now give a precise definition of universal convex bodies.

\A

\noindent {\bf Definition.} The convex body $K\in{\cal K}^n$ is
called {\em universal} if $\pi_m h_K \not= 0$ for all $m\in{\mathbb N}_0$.
The body $K\in{\cal K}^n$ is {\em centrally universal}\, if
$\pi_m h_K \not= 0$ holds for all even numbers $ m\in{\mathbb N}_0$.

\A

Since the space ${\cal H}^n_m$ and the scalar product on $C(S^{n-1})$ are
invariant under rotations, we have $\pi_m h_{\lambda\vartheta K} =
\lambda \vartheta (\pi_m h_K)$ for every $\vartheta \in SO(n)$ and
every $\lambda \geq 0$. Therefore, the property of being
universal (centrally universal) is invariant under rotations and
dilatations.

Let $K\in{\cal K}^n$ be a $k$-dimensional convex body, $k \geq 2$.
Because of the rotation invariance just mentioned and the
translation equivariance of the Steiner point map, it is no loss
of generality to assume that $K\subset {\mathbb
R}^k\subset{\mathbb R}^n$. If $K$ is universal in the sense of
the preceding definition, we say that $K$ is {\em universal in
${\mathbb R}^n$}. Alternatively, we can consider $K$ as a subset
of the Euclidean space ${\mathbb R}^k$. If it is universal there,
that is, in the sense of the definition with $n=k$, we say that
$K$ is {\em universal in ${\mathbb R}^k$}. The following was
proved in \cite[\S5]{Sch74}.

\A

\noindent {\bf Lemma 1.} {\em If $K$ is universal in ${\mathbb
R}^k$, then $K$ is universal in ${\mathbb R}^n$. }

\A

Before we turn to the proof of Theorem 2, we recall that the set
of convex bodies $L\in{\cal K}^n$ with $h_L \in \mathcal{H}^n$ is dense in
$\mathcal{K}^n$, see \cite[ p. 160]{Sch93}. Similarly, the set
of all $L\in{\cal K}_c^n$ with $h_L \in \mathcal{H}^n$ is dense in
$\mathcal{K}_c^n$, and the set of all $L\in{\cal K}_o^n$ with $h_L \in \mathcal{H}^n$ is dense in
$\mathcal{K}_o^n$.

\A

\noindent {\em Proof of Theorem} 2. (a) Suppose $K$ is not universal, i.e.,
$\pi_m h_K=0$ for some $m$. If $m=0$, then $K$ is one-pointed,
hence we can assume that $m \ge 1$. Let $Y_m$ be a non-zero
spherical harmonic of order $m$. There is a constant $c > 0$ such
that $c + Y_m = h_M$ for some convex body $M$, see \cite[Lemma 1.7.9]{Sch93}.
Every function $f$ in the closure of the vector space
spanned by the functions $h_{\vartheta K}$, $\vartheta \in SO(n)$,
satisfies $\pi_m f =0$. Since this does not hold for $h_M$, the
set of generalized ${\cal M}$-bodies is not dense in ${\cal K}^n$.

Now let $K$ be universal. By the preceding remark, it is enough
to show that for every convex body $L$ with $h_L \in
\mathcal{H}^n$ there are $\mathcal{M}$-bodies $T_1, T_2$ such
that $L + T_1 = T_2$. This was proved in \cite{Sch96} for the
case where $K$ is a triangle with at least one of its angles an irrational multiple of $\pi$. For an arbitrary universal
convex body $K$, the proof is almost verbally the same. We
sketch the argument, for the reader's convenience and since we have to
explain the modifications in part (b). Let
\begin{equation}\label{2.1}
h_L=\sum \limits_{m=0}^k \sum \limits_{j=1}^{N(n,m)} a_{mj} Y_{mj}
\end{equation}
and define $c_{mj}=(h_K,Y_{mj}).$
Since $K$ is universal, for each $m \in \mathbb{N}_0$ there is an
index $j(m)$ such that $c_{mj(m)} \neq 0$. With the functions $t_{ij}^m$ from (\ref{tij}), we define
\begin{equation}\label{2.2a}
g := N(n,m) \sum \limits_{m=0}^k \sum \limits_{i,j=1}^{N(n,m)}b_{ij}^mt_{ij}^m ,
\end{equation}
where $b_{ij}^m = a_{mi}c_{mj(m)}^{-1}$ for $j=j(m)$ and 0
otherwise. By (\ref{form1}), we have
\begin{eqnarray*}
\int_{SO(n)} h_{\vartheta K}(u)g(\vartheta)\,\D \nu(\vartheta) & = &
N(n,m) \sum \limits_{m=0}^k \sum \limits_{i,j=1}^{N(n,m)}
\int_{SO(n)} h_{\vartheta K}(u)b_{ij}^mt_{ij}^m(\vartheta)\,\D \nu(\vartheta)\\
 & = & \sum \limits_{m=0}^k \sum \limits_{i,j=1}^{N(n,m)}
 b_{ij}^m (h_K,Y_{mj})Y_{mi}(u)\\
 & = & \sum \limits_{m=0}^k \sum \limits_{i=1}^{N(n,m)}
 a_{mi}Y_{mi}(u)=h_L(u).
\end{eqnarray*}
Splitting $g$ into positive and negative parts, we obtain that $L$ is a generalized ${\cal M}$-body.

(b) If $K$ is not centrally universal, then $\pi_m h_K$ for some even number $m$.
With $M$ constructed as in part (a), we have $M\in{\cal K}_o^n$, but $M$
cannot be approximated by bodies from ${\cal M}$, since all bodies $L\in {\cal M}$
satisfy $\pi_mh_L=0$.

Suppose that $K$ is centrally universal. First let $K$ have its centre different from $0$. Then $\pi_1h_K\not=0$.
Let $L\in{\cal K}_c^n$ and $h_L\in{\cal H}^n$. Let $h_L$ be represented by (\ref{2.1}); then $a_{mj}=0$ for all odd $m\not=1$. For even numbers $m$ and for $m=1$, we can define
$c_{mj}$ and $b^m_{ij}$ as in part (a), as well as the function $g$, where now $b^m_{ij}:=0$ for odd $m\not=1$.

If $K$ has centre $0$, then $\pi_1h_K=0$ and we choose also $b^m_{ij}:=0$ for $m=1$. The proof can now be completed as in part (a).
\hfill$\Box$

\A

\section{Linear transformations}

In this section, we investigate the behaviour of $\pi_m h_{AK}$ under linear transformations $A$.

\A

\noindent{\bf Lemma 2.} {\em Let $f,g:{\mathbb
R}^n\setminus\{0\}\to{\mathbb R}$ be continuous functions, let
$f$ be positively homogeneous of degree $1$ and $g$ positively
homogeneous of degree $-(n+1)$. Then, for every $A\in GL(n)$,}
\[ \int_{S^{n-1}} f(Av)g(v)\,\D\sigma(v) =\frac{1}{|\det A|}
\int_{S^{n-1}} f(v)g(A^{-1}v)\,\D\sigma(v).\]

\A

\noindent{\em Proof.} We can take advantage of the known
transformation behaviour of the dual mixed volume
$\tilde{V}_{-1}(K,L)$. It is defined for two star bodies
$K,L\subset{\mathbb R}^n$ (i.e., compact sets, starshaped with
respect to the origin and with continuous radial functions) by
\[ \tilde{V}_{-1}(K,L) = \frac{1}{n} \int_{S^{n-1}} \rho(K,u)^{n+1}\rho(L,u)^{-1}\,\D\sigma(u),\]
where $\rho(K,u)=\max\{\lambda \geq 0: \lambda u \in K\}$ denotes
the radial function of $K$. Decomposing $f$ into its positive and
negative part and adding, say, the function $x\mapsto\|x\|$ to
both components, we can write $f= f^+ - f^-$, where the functions
$f^+,f^-$ are positive, continuous, and positively homogeneous of
degree $1$, on ${\mathbb R}^n\setminus\{0\}$. Let $L_+, L_-$ be
the star bodies with radial functions
$\rho(L_+,\cdot)=(f^+)^{-1}$ and $\rho(L_-,\cdot)=(f^-)^{-1}$,
then
\[ f = \rho(L_+,\cdot)^{-1} - \rho(L_-,\cdot)^{-1}.\]
Similarly, there are star bodies $K_+,K_-$ with
\[ g = \rho(K_+,\cdot)^{n+1} - \rho(K_-,\cdot)^{n+1}.\]
The assertion of the lemma now follows from the fact that $\rho(K,Av)=\rho(A^{-1}K,v)$ and that
\[ \tilde{V}_{-1}(AK,AL) = |\det A| \tilde{V}_{-1}(K,L), \]
which is proved in Lutwak \cite[Lemma (7.9)]{Lut90}.\hfill$\Box$

Lemma 2 is applied in the following argument. For $m\in{\mathbb
N}_0$ and $j\in\{1,\dots,N(n,m)\}$ we define
\[ \check{Y}_{mj}(x):= \frac{1}{\|x\|^{n+1}}Y_{mj}\left(\frac{x}{\|x\|}\right) \qquad \mbox{for }x\in{\mathbb R}^n\setminus\{0\}.\]
Let $K\in{\cal K}^n$ be fixed, and let $A\in GL(n)$. Since
$h(AK,v)=h(K,A^{T}v)$, Lemma 2 yields
\[ (h_{AK},Y_{mj}) =\frac{1}{|\det A|} \int_{S^{n-1}} h(K,v)\check{Y}_{mj}(A^{-T}v)\,\D\sigma(v).\]
Since the spherical harmonic $Y_{mj}$ is the restriction of a
polynomial on ${\mathbb R}^n$ to the sphere $S^{n-1}$, the
function
\[ A \mapsto\frac{1}{|\det A|}\check{Y}_{mj}(A^{-T}v) \]
is real analytic on the connected component of the identity.
(We identify $GL(n)$, via matrix description, with an open
subset of $\R^{n^2}$.)
The convergence of its power series is uniform on every compact
subset and uniform in $v$. Using the compactness of $S^{n-1}$, we
see that the function defined by
\[ A \mapsto (h_{AK},Y_{mj}), \]
(for given $K$) is real analytic. We will mostly apply this with
linear maps $A(\lambda)$ which, with respect to the standard
orthonormal basis of ${\mathbb R}^n$, have diagonal matrices ${\rm
diag}(1,\lambda,\dots,\lambda)$ or (in Section 4) ${\rm
diag}(1,1,\lambda,\dots,\lambda)$, with $\lambda\in I$. Here $I$
is any open interval $(0,a)$ with $a>1$. Then the function defined
by
\begin{equation} \label{fmj}
f_{mj}(\lambda) := (h_{A(\lambda)K},Y_{mj}), \qquad \lambda\in I,
\end{equation}
is real analytic on $I$.

We are now in a position to prove part (b) of Theorem 3.
This case exhibits already the basic idea for the proof of the general result.

\A

\noindent{\em Proof of Theorem} 3 (b). Let $K\in{\cal K}_c^n$ be a non-trivial convex body.
For $\lambda\in I$, let $A(\lambda)\in GL(n)$ be defined by
$A(\lambda):(x_1,\dots,x_n)\mapsto(x_1,\lambda x_2,\dots,\lambda
x_n)$ (where $x_1,\dots,x_n$ are coordinates with respect to the
standard orthonormal basis of ${\mathbb R}^n$). We may assume
that the orthogonal projection of $K$ to the first coordinate
axis is a segment $S$ of positive length. We have
$\lim_{\lambda\to 0} A(\lambda)K=S$ in the Hausdorff metric and
hence
\[ \lim_{\lambda\to 0} (h_{A(\lambda)K},Y_{mj}) = (h_S,Y_{mj})\]
for all $m,j$. Let $m$ be even. It is well known that, for all
$Y_m\in{\cal H}^n_m$, $u \in S^{n-1}$,
\[ \int_{S^{n-1}} |\langle u,v \rangle|Y_m(v)\,\D\sigma(v) = a_m Y_m(u)\]
with $a_m\not=0$ (see, e.g., \cite[ p. 185]{Sch93}). Since
$h_S= |\langle x,\cdot \rangle|+\langle y,\cdot \rangle$ for suitable $x, y \in \mathbb{R}^n$ it
follows that $S$ is centrally universal. Therefore, for each even
$m\in{\mathbb N}$ there is at least one index
$j(m)\in\{1,\dots,N(n,m)\}$ for which $(h_S,Y_{mj(m)}) \not=0$.
This implies that the function $f_{mj(m)}$, defined by
(\ref{fmj}), does not vanish identically on $I$. Since it is real
analytic, its zeros are isolated, thus there is an at most
countable subset $Z_m\subset I$ such that
$f_{mj(m)}(\lambda)\not=0$ for $\lambda\in I\setminus Z_m$. For
such $\lambda$, we have $\pi_m h_{A(\lambda)K}\not=0$ (this holds
trivially for $m=0$ and all $\lambda\in I$, since $A(\lambda)K$
has positive mean width). If now $U\subset {\mathbb R}$ is a
given neighbourhood of $1$, there exists a number \[\lambda\in
U\setminus\bigcup_{m\in{\mathbb N},\, 2|m} Z_m.\] It satisfies
$\pi_m h_{A(\lambda)K}\not=0$ for all even $m$, hence
$A(\lambda)K$ is centrally universal. Such a map $A(\lambda)\in
GL(n)$ can be found in any prescribed neighbourhood of the
identity. This completes the proof of Theorem 3 (b). \hfill$\Box$

\section{Proof of Theorem 3 (a)}

For the proof of Theorem 3 (a), we first
treat the case $n=2$. Let $K\in{\cal K}^2$ be a non-symmetric convex body; then $K$ has interior points.

As usual, we parameterize $S^1$ by an angle, writing
$u_{\varphi}:= (\cos \varphi,\sin\varphi)$, and by slight abuse
of notation, for a function $f$ on $S^1$ we write
$f(u_\varphi)=f(\varphi)$. The space ${\cal H}^2_m$ is spanned by
the functions $\cos m\varphi$ and $\sin m\varphi$, thus, in
complex notation,
\[ \pi_m h_K =0 \enspace\Leftrightarrow \enspace\int_0^{2\pi} h(K,\varphi)\,\E^{\I m\varphi} \,\D\varphi =0.\]
For $m \in {\mathbb N}$, we define a map $F_m(K,\cdot): GL(2)^+\to
{\mathbb C}$ by
\[ F_m(K,A) := \int_0^{2\pi} h(AK,\varphi)\,\E^{\I m\varphi}\,\D\varphi \qquad\mbox{for }A\in GL(2)^+,\]
where $GL(2)^+$ denotes the connected component of the identity in $GL(2)$.
As explained in Section 3, the function $F_m(K,\cdot)$ (interpreted as
a function on an open subset of $\R^4$) is real analytic. We make use of
the following fact (for a proof, see \cite[Lemma 5]{BG06}).

\A

\noindent {\bf Lemma 3.} {\em Let $f: U \rightarrow
\mathbb{R}$ be a real analytic function on an open subset of $\R^k$. Then the zero set of $f$
has Lebesgue measure zero, unless $f$ vanishes identically.}

\A

In view of Lemma 3, the following auxiliary result will be crucial.

\A

\noindent{\bf Lemma 4.} {\em The relation
\[ F_m(K,\cdot) \equiv 0 \]
does not hold for any odd integer $m\ge 1$.}

\A

\noindent For the proof, we assume that the assertion were false.
Then there exists a smallest odd integer $m\ge 1$ with
$F_m(K,\cdot) \equiv 0$.

\noindent{\em First Case:} $m\ge 5$. By the choice of $m$, we have
$F_{m-2}(K,\cdot)\not \equiv0$, hence there exists a map $A_0\in
GL(2)^+$ with $F_{m-2}(K,A_0)\not=0$, equivalently
$F_{m-2}(A_0K,{\rm Id})\not=0$. We still have $F_m(A_0K,\cdot)
\equiv0$. We may now replace $K$ by $A_0K$ and change the notation.
Thus, we can assume that
\begin{equation}\label{3.1}
F_m(K,\cdot) \equiv0,\qquad F_{m-2}(K,{\rm Id})\not=0.
\end{equation}
\noindent{\em Second Case:} $m =1$ or $m=3$. Since $K$ is not
symmetric, we cannot have $F_k(K,{\rm Id})=0$ for all odd
integers $k\ge 3$. A fortiori, $F_k(K,\cdot) \equiv 0$ cannot hold for
all odd integers $k\ge 3$. Therefore, there exists an odd integer
$k\ge 3$ such that $F_{k-2}(K,\cdot) \equiv 0$, but
$F_k(K,\cdot)\not\equiv 0$. As in the first case, we can replace
$K$ by $A_0K$, with a linear map $A_0$, so that after a change of notation we have
\begin{equation}\label{3.2}
F_{k-2}(K,\cdot) \equiv 0, \quad\mbox{but} \quad F_k(K,{\rm Id})\not=0.
\end{equation}

For a while, both cases are now treated together. To $K$ we will
apply rotations $R(\alpha)$ and linear maps $A(\lambda)$, given
by matrices
\[ \left(\begin{array}{rr} \cos\alpha & \sin\alpha\\ -\sin\alpha & \cos\alpha\end{array}\right) \qquad \mbox{and} \qquad \left(\begin{array}{rr} 1 & 0\\ 0 & \lambda\end{array}\right),\]
respectively, with $\alpha$ in an open neighbourhood $U$ of $0$
and $\lambda\in I$. We have
\[ h(A(\lambda)K,u_{\varphi}) = h(K,A(\lambda)u_{\varphi}) = \|A(\lambda)u_{\varphi}\|h(K,u_{\psi}) \]
with
\[ u_{\psi}:=\frac{A(\lambda)u_{\varphi}}{\|A(\lambda)u_{\varphi}\|}=\frac{(\cos\varphi,\lambda\sin\varphi)}{\sqrt {\cos^2\varphi +\lambda^2\sin^2\varphi}}.\]
Then
\[ u_{\varphi} = \frac{(\lambda\cos \psi,\sin\psi)}{\sqrt{\lambda^2\cos^2\psi+\sin^2\psi}},\]
thus
\[ \E^{\I m\varphi} = \frac{(\lambda\cos\psi+\I\sin\psi)^m}{(\lambda^2\cos^2\psi+\sin^2\psi)^{\frac{m}{2}}}\]
and
\[ \frac{\D\varphi}{\D\psi} = \frac{\lambda}{\lambda^2\cos^2\psi+\sin^2\psi}.\]
Substituting $\varphi$ by $\psi$ in the integral defining
$F_m(K,A)$, we get
\[ F_m(K,A(\lambda)) = \lambda^2\int_0^{2\pi} h(K,\psi)\frac{(\lambda\cos\psi+\I\sin\psi)^m}{(\lambda^2\cos^2\psi+\sin^2\psi)^{\frac{m+3}{2}}}\,\D\psi.
\]
By (\ref{3.1}), this integral vanishes for all $\lambda\in I$.
Therefore, also its derivatives with respect to $\lambda$ vanish.
For
\[ f(\lambda)= \frac{(\lambda\cos\psi+\I\sin\psi)^m}{(\lambda^2\cos^2\psi+\sin^2\psi)^{\frac{m+3}{2}}},\]
we obtain
 \[ 0=-f'(1) = \frac{3}{2}\,\E^{\I m\psi} +\frac{3-m}{4}\,\E^{\I (m-2)\psi} +\frac{3+m}{4}\,\E^{\I (m+2)\psi}.\]
Since $F_m(K,{\rm Id})=0$, this yields
\begin{equation}\label{3.3}
\int_0^{2\pi} h(K,\psi)[(3-m)\,\E^{\I (m-2)\psi}+(3+m)\,\E^{\I (m+2)\psi}]\,\D\psi=0.
\end{equation}
According to (\ref{3.1}), we also have $F_m(R(\alpha)K,\cdot)
\equiv 0$ for each angle $\alpha\in U$, hence (\ref{3.3}) holds
if $K$ is replaced by $R(\alpha)K$. Since
$h(R(\alpha)K,\psi)=h(K,\psi -\alpha)$, we see after a
substitution that (\ref{3.3}) holds with $\psi$ in the
exponentials replaced by $\psi+\alpha$. Since the functions
$\E^{\I (m-2)\alpha}$ and $\E^{\I (m+2)\alpha}$ are linearly
independent on $U$, we deduce that
\begin{equation}\label{3.4}
F_{m-2}(K,{\rm Id}) =0 \quad\mbox{if }m\not= 3,\qquad F_{m+2}(K,{\rm Id}) =0.
\end{equation}

Now we distinguish between the two cases considered above. If
$m\ge 5$, then the first relation of (\ref{3.4}) yields
$F_{m-2}(K,{\rm Id}) =0 $, which contradicts the second relation
of (\ref{3.1}). If $m=1$ or $m=3$, then we note that the first relation of
(\ref{3.2}) gives $F_{k-2}(K,\cdot)\equiv 0$. Therefore, the
second relation of (\ref{3.4}) holds also with $m$ replaced by
$k-2$, but this contradicts the second relation of (\ref{3.2}).
This completes the proof of Lemma 4. \hfill$\Box$

Let $m\in{\mathbb N}_0$ be an integer. If $m$ is odd, it follows from Lemma 4
that the real analytic function $F_m(K,\cdot)$ does not vanish
identically. If $m$ is even, the same result follows as in the
proof of Theorem 3 (b). Hence, the set of zeros of $F_m(K,\cdot)$ has
Lebesgue measure zero, by Lemma 3. Therefore, in any given neighbourhood of the
identity in $GL(2)$, we can find a linear map $A$ with
$F_m(K,A)\not=0$ for all $m\in{\mathbb N}_0$. The convex body $AK$
is universal. This completes the proof of Theorem 3 for $n=2$ and
non-symmetric convex bodies.

Finally, we assume that $n\ge 3$ and that $K\in {\cal K}^n$ is a
non-symmetric convex body. Then $K$ has dimension at least two.
There exists (see, e.g., Gardner \cite[Corollary 3.1.5]{Gar95}) a
two-dimensional subspace, without loss of generality the space
${\mathbb R}^2 \subset {\mathbb R}^n$, such that the orthogonal
projection $K'$ of $K$ to ${\mathbb R}^2$ is non-symmetric. In a
given neighbourhood of the identity of $GL(n)$ we can find an
affine transformation $B$ which maps ${\mathbb R}^2$ into itself,
leaves the orthogonal complement of ${\mathbb R}^2$ in ${\mathbb
R}^n$ pointwise fixed, and is such that $BK'$ is universal in
${\mathbb R}^2$. By Lemma 1, $BK'$ is universal in ${\mathbb
R}^n$. Moreover, $BK'$ is the image of $BK$ under the orthogonal
projection to ${\mathbb R}^2$. Assuming that ${\mathbb R}^2$ is
spanned by the first two vectors of the standard orthonormal
basis of ${\mathbb R}^n$, we define linear maps $A(\lambda)$ by
$A(\lambda):(x_1,\dots,x_n)\mapsto(x_1,x_2,\lambda
x_3,\dots,\lambda x_n)$, $\lambda\in I$. As in the proof of
Theorem 3 (b), we have
\[ \lim_{\lambda\to 0} (h_{A(\lambda)BK},Y_{mj}) = (h_{BK'},Y_{mj}).\]
Since $BK'$ is universal in ${\mathbb R}^n$, for each $m$ there
exists an index $j(m)$ with $(h_{BK'},Y_{mj(m)}) \not=0$. Thus,
the function $\lambda \mapsto (h_{A(\lambda)BK},Y_{mj(m)})$,
$\lambda\in I$, does not vanish identically. The proof can now be
completed as it was done in Section 3 for centrally symmetric
bodies. \hfill$\Box$

\noindent Professor R. Schneider,\\Mathematisches Institut,\\Albert-Ludwigs-Universit\"at,\\
Eckerstra{\ss}e 1,\\D-79104 Freiburg i. Br.,\\Germany.\\E-mail: rolf.schneider@math.uni-freiburg.de

\vspace{4mm}

\noindent Dr. F. E. Schuster,\\Institut f\"ur Diskrete Mathematik
und Geometrie,\\Technische Universit\"at Wien,\\Wiedner
Hauptstra{\ss}e 8/104,\\A-1040 Wien,\\Austria.\\E-mail:
fschuster@osiris.tuwien.ac.at
\end{document}